\documentclass[12pt]{amsart}

\usepackage{graphicx,psfrag}

\newtheorem{prop}{Proposition}
\newtheorem{defin}{Definition}
\newtheorem{thm}{Theorem}

\newtheorem{corol}{Corollary}\newtheorem{cor}{Corollary}
\newtheorem{lem}{Lemma}

\newtheorem{rem}{Remark}

\def\qed{\hfill\vrule height 5pt width 5 pt depth 0pt}

\renewcommand{\a}{{\mathfrak a}}\newcommand{\A}{{\mathcal A}} 

\renewcommand{\b}{{\mathfrak b}}

\newcommand{\C}{{\mathbb C}}

\newcommand{\f}{{\mathfrak f}}

\newcommand{\g}{{\mathfrak g}}

\newcommand{\HH}{{\mathbb H}}

\renewcommand{\k}{{\mathfrak k}}

\newcommand{\n}{{\mathfrak n}}\newcommand{\N}{{\mathbb N}}

\newcommand{\p}{{\mathfrak p}}

\newcommand{\R}{{\mathbb R}}

\newcommand{\ttt}{{\mathcal T}}

\newcommand{\Z}{{\mathbb Z}}

\newcommand{\adj}{\text{Ad}}

\newcommand{\diam}{\text{diameter}}

\newcommand{\grexp}{\text{exp}}
\newcommand{\riexp}{\text{Exp}}
\newcommand{\Eig}{\text{Eig}}

\newcommand{\iso}{\text{Iso}}   

\newcommand{\mac}{\text{SL}(2,\R)}   
\newcommand{\mcz}{\text{SL}(2,\Z)}   

\newcommand{\rk}{\text{rk}} 

\newcommand{\tM}{\tilde{M}}
\newcommand{\tF}{\tilde{F}}
\newcommand{\tG}{\tilde{G}}
\newcommand{\tX}{\tilde{X}}

\newcommand{\tp}{\tilde{p}}

\newcommand{\tx}{\tilde{x}}
\newcommand{\ty}{\tilde{y}}
\newcommand{\tz}{\tilde{z}}
\newcommand{\tw}{\tilde{w}}

\newcommand{\tga}{\tilde{\ga}}

\newcommand{\ep}{\varepsilon}
\newcommand{\la}{\lambda}
\newcommand{\ga}{\gamma}\newcommand{\Ga}{\Gamma}

\begin{document}

\bibliographystyle{plain}

\title[Symmetric Spaces and Security]
{Connecting geodesics and security of configurations in compact locally symmetric spaces}

\author{Eugene Gutkin and Viktor Schroeder}

\address{University of California and IMPA\\ 
Estrada Dona Castorina 110\\
Rio de Janeiro/Brasil 22460-320\\
Department of Mathematics\\  
University of Z\"urich\\
CH-8057 Z\"urich\\
Switzerland}

\email{gutkin@impa.br,vschroed@math.unizh.ch}

\date{\today}

\begin{abstract} 
A pair of points in a riemannian manifold makes a secure configuration
if the totality of geodesics connecting them can be blocked by a finite set.
The manifold is secure if every configuration is secure.
We investigate the security of compact, locally symmetric spaces.    
\end{abstract}

\maketitle

\section{Introduction: The setting and the main results}       \label{intro}
Let $M$ be a complete riemannian manifold.\footnote{
Most of the preliminary material is valid in greater generality \cite{Gut03,Gut04}.
Since  locally symmetric spaces fit into the riemannian framework, we will restrict
our discussion to this setting.} 
By a geodesic $\ga\subset M$ we will mean a geodesic curve $t\mapsto\ga(t)$, where
$t\in I$ is the arclength parameter, and $I\subset\R$ is an arbitrary interval. 
Mostly, we will be concerned with the situation $I=[a,b]$, i. e., $\gamma$ is a 
geodesic segment with the endpoints $x=\ga(a),y=\ga(b)$.
\begin{defin}  \label{connect_def} 
Let $x,y\in M$ be arbitrary points. 
A {\em connecting geodesic} is a geodesic segment $\gamma$
with the endpoints $x,y$, and such that 
$\gamma$ does not contain either $x$ or $y$ in its interior.
\end{defin}
By a {\em configuration} in $M$ we will mean any unordered pair
of points, $\{x,y\}$. Let $\sigma:M\times M\to M\times M$ be the involution
$\sigma(x,y)=(y,x)$. The {\em space of configurations} is the quotient
$C(M)=(M\times M)/\sigma$. Thus, $C(M)$ is the symmetric square of $M$,
and it inherits from $M$ a topology, a differentiable structure, a measure class, etc.
If $M$ carries a group action, then the group, $G$, naturally acts on $C(M)$.
We will say that two configurations $\{x,y\},\{x',y'\}\in C(M)$
are {\em conjugate} if $\{x',y'\}=g\cdot\{x,y\}$ for some $g\in G$.

Let $z\in M$ and let $\ga$ be any geodesic. We say that $\ga$ {\em passes through} $z$ 
if $\ga$ contains $z$ in its interior.
Let $\Ga$ be any collection of geodesics in $M$ and let $F\subset M$ be a subset.
We say that $F$ is a {\em blocking set} for $\Ga$ if every geodesic in $\Ga$ 
passes through a point of $F$.

\begin{defin}  \label{secur_def} 
Let $\Gamma(x,y)$ be the collection of connecting geodesics for a configuration
$\{x,y\}$. We say that $\{x,y\}$ is a 
{\em secure configuration} if there exists a finite blocking set for $\Gamma(x,y)$.
Otherwise the configuration is {\em insecure}. The manifold $M$ is secure if
every configuration $\{x,y\}$ is secure.
\end{defin}

If $M$ is secure, and any collection $\Gamma(x,y)$ can be blocked by
a set of at most $n$ points, we say that $M$ is {\em uniformly secure}.
The smallest such $n$ is the {\em security threshold} of $M$.

In a geometric optics interpretation, a configuration is secure
if one of the points can be shaded from the light emanating from the other
by a finite number of point screens. Thus, our setting is closely related to the 
problem of illumination \cite{Rau}.
Another obvious interpretation of Definition~\ref{secur_def} suggested the name ``security".

Since the security of configurations concerns the global properties of geodesics,
it is instructive to investigate the possibilities and to compare 
various spaces from this viewpoint. The work \cite{Gut04} did this for a particular class
of planar polygons, the {\em lattice polygons}. The geodesics (i.e., the {\em billiard
orbits}) in a lattice polygon have a striking behavior: a geodesic is either
finite or uniformly distributed. It is the direction of the geodesic that
determines which of the two possibilities happens.

The name {\em lattice polygons} is due to the fact that any polygon, $P$, in  this
class defines a nonuniform lattice $G(P)\subset\mac$. Such a lattice is either
{\em arithmetic} (i.e., commensurable with $\mcz\subset\mac$) or {\em nonarithmetic}.
Accordingly, $P$ is either an arithmetic or a nonarithmetic polygon \cite{GJ}.
By a theorem in \cite{Gut04}, a lattice polygon is secure iff it is arithmetic.
Regular polygons are lattice polygons \cite{Ve1}. By \cite{Gut04} and \cite{Ve1},
a regular $n$-gon is secure iff $n=3,4,6$. Thus, any regular $n$-gon other than
the equilateral triangle, the square and the regular hexagon, is insecure.

A space $M$ is insecure iff $M$ has at least one insecure configuration.
It is natural to analyze insecure spaces by classifying their configurations
from the security viewpoint; the paper \cite{Gut05} does this for nonarithmetic 
lattice polygons of small genus. In view if the results of \cite{Gut05},
it is plausible that almost all configurations in a nonarithmetic lattice polygon are insecure.

In this work we investigate the security of a well known class of riemannian manifolds:
compact, locally symmetric spaces. Let $M$ be one. Then $M=S/\Ga$, where
$S$ is a simply connected symmetric space, and $\Ga$ is a discrete, cocompact group
of isometries freely acting on $S$. The space $S$  uniquely decomposes, 
$S=S_0\times S_-\times S_+$, into a product of 
simply connected symmetric spaces of {\em euclidean type, noncompact type, and 
compact type} respectively \cite{He1}. If $M=S/\Ga$ where $S$ belongs to one of
the three types, we say that $M$ is a compact, locally symmetric space of that type.

\medskip

We will now formulate the main results of this work.

\medskip


\noindent 1. Any configuration in a compact, locally symmetric space of the 
noncompact type is insecure. See Theorem~\ref{nct_thm}.

\noindent 2. Let $M$ be a (necessarily compact) locally symmetric space 
of compact type. We define the notion of {\em regular/singular configurations}.
The set of regular configurations is open and dense. Then:
i) Any regular configuration is secure. The security threshold 
of regular configurations is $2^{\rk(M)}$;
ii) There are always singular configurations which are insecure. 
See Theorem~\ref{block_thm}.

\noindent 3. Let $M$ be an arbitrary compact, locally symmetric space.
Then $M$ is secure iff it is of euclidean type. If $M$ is of euclidean type, then
it is uniformly secure, and its security threshold is bounded in terms of
$\dim(M)$. See  Theorem~\ref{gen_sec_thm} and Corollary~\ref{gen_sec_cor}.

\medskip

The organization of the paper is as follows. In Section~\ref{set} we collect  basic
facts, in particular on the security and coverings. There we also establish the
security of locally symmetric spaces of euclidean type. 
In Section~\ref{nct} we investigate the configurations
in a compact, locally symmetric space of the noncompact type, and prove Theorem~\ref{nct_thm}.
In Section~\ref{comp_type} we study the security of spaces with a compact
group of isometries. Theorem~\ref{act_thm} gives a sufficient condition for
insecurity of such spaces. Then we apply this material to 
locally symmetric space of the compact type, 
and prove Theorem~\ref{block_thm}. Subsections~\ref{exam1_sub} and~\ref{exam2_sub}
illustrate   Theorem~\ref{block_thm} by examples. In section~\ref{exam1_sub}
we consider compact symmetric spaces of rank one.
Using Theorem~\ref{block_thm}, we characterize their configurations from the
security viewpoint. In section~\ref{exam2_sub} we investigate
compact, semisimple Lie groups endowed with double-invariant
riemannian metrics (symmetric spaces of {\em type $II$} 
\cite{He2}). Theorem~\ref{secur_group_thm} is a direct corollary of Theorem~\ref{block_thm}.
In Section~\ref{gener_loc_symm} we consider arbitrary compact, 
locally symmetric spaces,  and prove Theorem~\ref{gen_sec_thm} and Corollary~\ref{gen_sec_cor}.

\section{Preliminaries}       \label{set}
We first discuss riemannian coverings from the security viewpoint.
This material is used in Section~\ref{comp_type}.
Then we recall the basic material on symmetric spaces, in general.
We will give more detailed presentations separately  
for the noncompact and the compact type.
See sections~\ref{notcomp_set} and~\ref{comp_set} respectively.

\subsection{Security and coverings}  \label{cov_sub}
Let $X,Y$ be complete, connected riemannian manifolds, and let $p:X\to Y$
be a differentiable mapping which is onto, and is a local isometry. 
Then $p$ is a topological covering. 
To limit our discussion to the security context, we assume that $X$ is compact,
and hence $p$ is a finite covering.\footnote{
The material below is valid in the more general framework of {\em geodesic coverings}
\cite{Gut03,Gut04}. Topological coverings suffice for our purpose, and
we restrict the discussion to them.}
\begin{prop} \label{cov_sec_prop}
Let $p:X\to Y$ be as above. A configuration $\{x,y\}$ in $Y$ is secure iff
all configurations $\{\tx,\ty\}$ in $X$, with $\tx\in p^{-1}(x),\,\ty\in p^{-1}(y),$
are secure.
\end{prop}
\begin{proof} Let $d$ be the degree of the covering, and let $x,y \in Y$ be arbitrary. Then
\begin{equation} \label{preim_eq}
p^{-1}(\Ga(x,y))=\cup_{\tx\in p^{-1}(x),\,\ty\in p^{-1}(y)}\Ga(\tx,\ty).
\end{equation}
Let $F\subset Y$ be a blocking set for $\{x,y\}$. Then, by eq.~\ref{preim_eq}, $\tF=p^{-1}(F)$
blocks the union of  $\Ga(\tx,\ty)$ over $\tx\in p^{-1}(x),\,\ty\in p^{-1}(y)$.  
Since $|\tF|=d|F|<\infty$, all these configurations are secure.
This proves one implication.
To prove the converse, set $p^{-1}(x)=\{\tx_1,\dots,\tx_d\},\,p^{-1}(y)=\{\ty_1,\dots,\ty_d\}$.
For $1\le i,j \le d$ let $\tF_{i,j}\subset X$ be a blocking set for $\Ga(\tx_i,\ty_j)$.
By eq.~\ref{preim_eq}, $F=p(\cup_{i,j}\tF_{i,j})$ blocks $\Ga(x,y)$, and 
$|F|\le\sum_{i,j}|\tF_{i,j}|<\infty$. \end{proof}

\medskip

The statement below is immediate from Proposition~\ref{cov_sec_prop} and its proof. 
\begin{corol} \label{cov_sec_cor}
Let $p:X\to Y$ be a covering of compact, riemannian manifolds. Then

\noindent 1. One of the manifolds is (uniformly) secure iff the other one is;

\noindent 2. Suppose that $Y$ is insecure, and let $\{x,y\}\in C(Y)$ be 
an insecure configuration. Then there exist $\tx\in p^{-1}(x),\,\ty\in p^{-1}(y)$
such that the configuration $\{\tx,\ty\}\in C(X)$ is insecure.
\end{corol}

\medskip

\subsection{Symmetric and locally symmetric spaces}  \label{basic_sub}
We will denote Lie groups by capital latin letters, and their Lie algebras
by the corresponding lower case gothic letters. Thus, if $G$ is a Lie group,
then $\g$ is the Lie algebra of  $G$. We denote by $G_0\subset G$
the connected component of identity. We refer the reader to \cite{He1,He2}
for the background on symmetric spaces, Lie groups, and Lie algebras.
See also \cite{Lo}.

A (riemannian, globally) symmetric space is a complete,
homogeneous riemannian 
manifold, $S=G/K$, where $G$ is a connected Lie group with an involutive
automorphism, $\sigma:G\to G$, and $K\subset G$ is (essentially) the
fixed point set of $\sigma$. We will use the same notation for the
induced automorphism of the Lie algebra.
The  automorphism, $\sigma:\g\to \g$, has eigenvalues $\pm 1$, and 
let $\g=\k + \p$ be the decomposition
of the Lie algebra into the eigenspaces of $\sigma$. The eigenspace $\p$
is naturally identified with the tangent space $T_oS$ at the reference
point $o=eK\in S$. The involution $\sigma:G\to G$ descends to an isometry
$s_o:S\to S$ such that $s_o^2=1$, $s_o(o)=o$ and $s_o|T_oS=-Id$.
Thus, $s_o$ is the geodesic symmetry of $S$ with respect to the reference point.
The action of $G$ on $S$ gives rise to the geodesic symmetries $s_x:S\to S$
where $x\in S$ is arbitrary.

The property of having a geodesic symmetry for every point can be used as 
a definition of symmetric spaces \cite{He2}. We assume
that $G$ acts faithfully on $S$, i. e., $G\subset\iso(S)$. Then $K$
is compact, and $G$ is a reductive Lie group. The Lie algebra $\g$
has a unique $\sigma$-invariant decomposition $\g=\g_0\oplus\g_-\oplus\g_+$
where $\g_0$ is the center of $\g$ and $\g_-,\,\g_+$ are noncompact and compact 
semisimple Lie algebras respectively. 
If $\g=\g_0$ ($\g=\g_-$, $\g=\g_+$), we say that the 
symmetric space $S$ is of euclidean type (noncompact type,  
compact type). A symmetric space of the euclidean type
satisfies $S_0=\R^n/\Ga$, where $\Ga\subset\R^n$ is a discrete subgroup.
The structure of a symmetric space of either type
is described via root decompositions of the corresponding Lie algebras.
We will recall this material separately for the spaces of noncompact
(section~\ref{notcomp_set}) and compact type (section~\ref{comp_set}).

An irreducible symmetric space necessarily belongs to one of the three types.
The general symmetric space  decomposes (at least locally) 
 as a cartesian product, $S=S_0\times S_-\times S_+$, of 
symmetric spaces of euclidean, noncompact, and 
compact type.\footnote{
This decomposition certainly exists (and is unique)
if $S$ is simply connected \cite{He1}.}
Locally symmetric (compact) spaces associated with the symmetric space $S=G/K$
are of the form $M=\Ga\setminus S$, where $\Ga\subset G$ is a discrete (cocompact) subgroup
freely acting on $S$. If $M$ is a locally symmetric space, and
$S$ belongs to a particular type, we will say that $M$ is a locally symmetric space
of the corresponding type. 

We dispose of the euclidean type in the subsection below. 
In the following two sections we study the security of 
locally symmetric spaces of the noncompact 
and the compact type respectively.

\subsection{Security of locally symmetric spaces of euclidean type}  
\label{euclid_sub}
A compact, locally symmetric space of euclidean type is of the form
$M^n=\R^n/\Ga$, where $\Ga\subset\iso(\R^n)$ is a cocompact,
freely acting, discrete subgroup.
A finite covering of $M^n$ is a flat torus of dimension $n$.

\begin{prop} \label{euclid_case}
Any compact, locally symmetric space $M^n$ of euclidean type is uniformly secure;
there is a bound on security thresholds of these spaces, depending only on $n$.
If $M^n$ is a flat torus, then the security threshold is $2^n$.
\end{prop} 
\begin{proof} The case $n=2$ is contained in \cite{Gut04}, Lemma 1. 
The same aproach works for any $n$. We outline it below.

By the Bieberbach theorem, $M^n$ has a finite covering by a flat torus;
moreover, the degree of the covering is bounded above in terms of $n$.
In view of Corollary~\ref{cov_sec_cor}, it suffices to consider the case
$M^n=\R^n/\Ga$, where $\Ga \subset \R^n$ is a lattice.

Affine transformations $M\to g\cdot M$ preserve the set of geodesics in $M$. Thus, the claim holds
for $M$ iff it holds for any $g\cdot M$. Using an appropriate  $g$, we can assume that
$M=\R^n/\Z^n$, the standard torus $T^n$ of $n$ dimensions. Let $o\in T^n$ be the origin.
By homogeneity, it suffices to consider the configurations $\{o,x\}$.

There is a one-to-one correspondence between the geodesics $\ga\in\Ga(o,x)$
and the straight segments $\tga_{x+z}$ in $\R^n$ connecting the origin $0\in\R^n$ with
the points $x+z,\,z\in\Z^n$. Let $\ga_{x+z}\in\Ga(o,x)$ be the corresponding
connecting geodesic. If $p:\R^n\to T^n$ is the projection, then $\ga_{x+z}=p(\tga_{x+z})$.
The midpoint of the segment  $\tga_{x+z}$  is $\frac{x}{2}+z/2\in\R^n$.
Set $\tF(x)=\{\frac{x}{2}+z/2:\,z\in\Z^n\}$.
Then the set $F(x)=p(\tF(x))\subset T^n$ is finite, and $|F(x)|=2^n$.
Thus, $2^n$ points suffice to block any $\Ga(o,x)$. On the other hand,
for a typical $x$, we cannot block $\Ga(o,x)$ with less than $2^n$ points.
We leave the verification of this to the reader. \end{proof}


\section{Compact locally symmetric spaces of noncompact type}    
\label{nct}
We begin by presenting preliminaries, and establishing notation.

\subsection{Symmetric spaces of noncompact type}       
\label{notcomp_set}
A symmetric space of noncompact type satisfies $S=G/K$, where
$G = \iso_0(S)$ is a noncompact, semisimple Lie group, and $K$ is a maximal compact subgroup.
The subgroup $K\subset G$ is defined up to conjugation; 
there is a one-to-one correspondence between the choices
of $K$ and the choices of a reference point in $S$.
It will be convenient to consider any point
$x \in S$ a reference point. 

We will denote by $\g = \k + \p$ the Cartan decomposition corresponding to
$x \in S$; here $\k$ is the Lie algebra of $K$, and $\p \simeq T_xS$.
The Riemannian exponential map, $\riexp:\p \to S$, 
and the Lie group exponential map, $\grexp:\g \to G$, are related by
$\riexp(H) = \grexp(H)\cdot x$. A {\em flat}, $X\subset S$, is a
totally geodesic submanifold, isometric to a euclidean space.
If $\a\subset \p$ is a maximal abelian subalgebra, 
then $\riexp(\a)=\grexp(\a)\cdot x\subset S$ is a maximal flat. Varying $x\in S$ and
$\a\subset \p$, we obtain all maximal flats in $S$.

We will use the standard material (and the standard notation) on
root systems and the root decompositions \cite{He1,Eb3}. 
Thus, a maximal abelian subalgebra $\a\subset\p$ gives rise to the root decomposition
\begin{equation}  \label{root_eq}
\g = \g_0 + \sum_{\la \in \Delta} \g_{\la}.
\end{equation} 
For a root $\la\in\Delta$, the root vector is the unique element
$H_{\la} \in \a$ such that $\la(H) = <H_{\la},H>$ for all $H\in \a$.
Weyl chambers are the connected components of
$(\a \setminus \bigcup _{\la \in \Delta} H_{\la}^{\perp})$.
Every vector $H \in \p$ is contained in a maximal abelian subalgebra
$\a\subset \p$, and all maximal abelian subalgebras of $\p$ are
$K$-conjugate. Thus, the set $\{\la(H):\la\in\Delta\}\subset\R$
does not depend on the choice of $\a$; it is the set of nontrivial
eigenvalues of the symmetric linear transformation
$ad(H):\g \to \g$ \cite[(2.7.1)]{Eb3}.
A vector $H\in\p$ is {\em regular} iff it is contained in a unique
maximal abelian subalgebra $\a$, iff $\la(H)\ne 0$ for all $\la\in\Delta$,
iff $H$ belongs to a unique Weyl chamber \cite{He1,Eb3}.

A symmetric space of noncompact type is a Hadamard
manifold, i. e., $S$ is a complete, simply connected riemannian manifold of
nonpositive sectional curvature.
We will use the standard properties of Hadamard manifolds \cite{Eb3}.
For two points $x,y \in S$ there exists a
unique geodesic, $[x,y]$, from $x$ to $y$. Let $J\subset\R$ be an
arbtrary interval, 
and let $g,h:J \to S$ be two geodesics parametrized by the arclength.
Then the distance function $t \mapsto d(g(t),h(t))$
on $J$ is convex. Two
geodesic rays $\ga , \sigma :\R_+ \to S$ are {\em asymptotic}
if $d(\ga(t),\sigma(t))$ is bounded as $t \to +\infty$.
The {\em ideal boundary}  $\partial_{\infty}S$ of the symmetric space $S$
is the set of classes of asymptotic geodesic rays.
We denote by $\ga(\infty) \in \partial_{\infty}S$ the boundary point
corresponding to the geodesic ray $\ga$.

Denote by $T^1S$ the unit tangent bundle of $S$,
and by $T^1_xS\subset T^1S$ its fiber at $x\in S$.
Let $v\in T^1_xS$. Let $H\in\p$ be the vector corresponding to $v$
via an isomorphism $\p\simeq T^1_xS$, and let $\a\subset\p$ 
be a maximal abelian subalgebra containing $H$. By preceding remarks,
the set $\Eig(v) = \{ \la(H) : \la \in \Delta \}$ is
well defined. This correpondence defines a continuous (with respect to
the Hausdorff distance) set valued function $\Eig(v)$ on $T^1S$.
Hence, the function
$\la_0^{+}(v)= \min \{ \mid l\mid : l \in \Eig(v) \}$ is well defined 
and continuous on $T^1S$. For any $v\in T^1S$ we have $\la_0^{+}(v)\ge 0$;
a vector $v$ is {\em regular} if $\la_0^{+}(v)> 0$. The regularity of elements
$v\in T^1S$ agrees with the regularity for vectors $H\in\p$.

The set-valued function $\Eig$ on $T^1S$ is invariant under the action of 
$G$ and under the geodesic flow. Therefore $\Eig(v),\,v\in T^1S$, is determined
by the geodesic $\ga=\riexp(tv)$. Moreover, $\Eig(v)$ depends only on
the point $\ga(\infty) \in \partial_{\infty}S$.
In view of these remarks, $\la_0^{+}$ uniquely descends
to a continuous function
$\la_0^{+}:\partial_{\infty}S\to\R_+\cup\{0\}$.
The notion of regularity for tangent vectors $v\in T^1_xS$
defines the regularity for boundary points $\xi \in \partial_{\infty}S$.
By remarks above, a point $\xi \in \partial_{\infty}S$ is
regular iff $ \la_0^{+} (\xi) > 0$.

A regular point $\xi \in \partial_{\infty}S$ and any point
$x \in S$ determine the {\em horocycle} $HC(\xi,x)$.
(Compare with \cite[pp.105-108]{Eb3}).
Let $\ga$ be the geodesic from $x$ to $\xi$, and let $v = \dot \ga(0)\in T^1_xS$.
Let $\g = \p + \k$ be the Cartan decomposition corresponding to 
$x$, and let $H \in \p$ be the vector corresponding to $v$. Then $H$ is regular,
and let $\a\subset\p$ be the unique maximal abelian subalgebra containing $H$.
Since $H$ is regular, it is contained in a unique Weyl chamber; let
$\Delta^{+} \subset \Delta$ be the corresponding set of positive roots, 
i. e. $\la \in \Delta^{+}$ iff $\la(H) > 0$. Set
$$
\n = \sum_{\la \in \Delta^{+}} \g_{\la},
$$
and let $N = \grexp(\n)\subset G$ be the corresponding nilpotent subgroup.
Then $HC(\xi,x) = N \cdot x\subset S$.

\begin{rem} \label{rem_horocycle}
{\em Let $x\in S$ be arbitrary.
Every horocycle containing $x$ has the form
$HC(k\xi,x)$, where $k\in K$, the isotropy group of $x$.
Therefore the set of horocycles passing  through $x$ is compact.}
\end{rem}

\vspace{3mm}

For $x\in S$ and $r>0$ let $B_r(x)\subset S$ be the ball of radius $r$ centered at $x$.
For $x,y\in S$ and $r>0$
set $SH(x,y,r)= \{u \in S : [x,u] \cap B_r(y) \neq \emptyset \}$.
The set $SH(x,y,r)\subset S$ is the {\em shadow} of the ball $B_r(y)$
produced  by the light emitted from $x$. Let $s>0$. We will call
the intersection  $SH(x,y,r) \cap B_s(y)$ a {\em restricted shadow}.


\subsection{Any configuration is insecure: Outline of the proof}
\label{noncomp_type}
 Let $M=  S/\Gamma$ be a compact, locally symmetric space of noncompact type.
Thus, $S=\tM$ is a simply connected, noncompact symmetric space,
and  $\Gamma \subset I(S)$ is the deck group of the covering
$\pi: S \to M$. The following is the main result of this section.

\begin{thm}    \label{nct_thm}
Let $M$ be a compact locally symmetric
space of noncompact type, let $x,y \in M$ and let
$F \subset M\setminus\{x,y\}$ be a finite set. Then there exists a geodesic
$h\in \Gamma(x,y)$ such that $h \cap F = \emptyset$.
\end{thm}

For the benefit of the reader, below we sketch a proof of Theorem~\ref{nct_thm}. 
First, we consider locally symmetric spaces of rank one, and
outline an argument that proves the claim in this case.
It is substantially simpler than the argument for  higher rank
locally symmetric spaces; it is especially transparent for
compact surfaces of constant negative curvature.
Then we outline the argument  for  higher rank locally symmetric spaces,
emphasizing the modifications and the difficulties that do not
arise in the rank one case.

\vspace{3mm}

Let $\rk (M) =1$. By \cite{BSW}, there exists a closed geodesic,
$\alpha \subset M$, such that $\alpha \cap (F \cup\{ x,y\})= \emptyset$.
Let $A \subset S$ be an infinite geodesic such that $\pi(A) = \alpha$.
For $\ep > 0$, we will denote by $T_{\ep}(X)$ the $\ep$-tube about $X$; for $X \subset M$
we set $\tX = \pi^{-1}(X) \subset S$.

Since $\alpha$ is compact, there exists $\ep> 0$ such that
$T_{\ep}(A) \cap \tF = \emptyset$. Let $\tp_i \in A,\,1\le i,$ 
be a sequence of points going to infinity.
Let $\tx \in \pi^{-1}(x)$, let $\tz_i=s_{\tp_i}(\tx)$, and set
$\gamma_i=[\tx,\tz_i]$.
Let $l_i=|\gamma_i|= d(\tx,\tz_i)$. We parametrize the geodesics $\gamma_i$
so that $\gamma_i(0)=\tx,\gamma_i(l_i)=\tz_i$.

By hyperbolicity, there exists $\rho > 0$  such that 
$\gamma_i((\rho,l_i-\rho))\subset T_{\ep}(A)$
for all $i$. Thus, the geodesics $\gamma_i$ are contained in the 
$\ep$-tube about $A$ except for, possibly, the first
and the last segments of length at most $\rho$, where $\rho$ does not depend
on $i$.
Let $\eta > 0$, let $\tw_i \in B_{\eta}(\tz_i)$ be arbitrary, and set $\gamma_{\tw_i}=[\tx,\tw_i]$.
Then the preceding claim holds for all geodesics $\gamma_{\tw_i}$
if $\eta > 0$ is sufficiently small. Note that $\gamma_i=\gamma_{\tz_i}$.

For $\eta > 0$ and $1\le i$, we set $SH_i(\eta)=SH(\tx,\tz_i,\eta)$
and $SH_i(\eta,r)=SH_i(\eta) \cap B_r(\tz_i)$.
Thus, the sets $SH_i(\eta)$ and $SH_i(\eta,r)$ are the shadow and
the restricted shadow respectively, in our context.

\begin{figure} \label{rankone}
\psfrag{x}{$\tilde{x}$}
\psfrag{p}{$\tilde{p}_{i}$}
\psfrag{z}{$\tilde{z}_{i}$}
\psfrag{y}{$\tilde{y}_{i}$}
\psfrag{A}{$A$}
\includegraphics{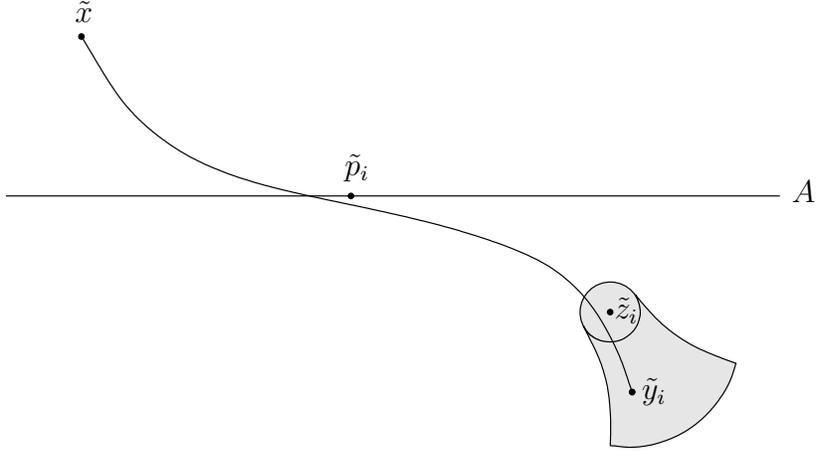}%
\caption{Construction of connecting geodesics in the rank one case}
\end{figure}

Since the curvature of $S$ is bounded away from
zero, the standard comparison arguments yield that for any $\eta>0$ 
there exists $r > 0$ such that each $SH_i(\eta,r)$ contains a ball in $S$
of radius $\diam (M)$.Thus, $\pi(SH_i(\eta,r))=M$. In particular, 
the set $SH_i(\eta,r)$ contains a point $\ty_i \in \pi^{-1}(y)$.
See figure~\ref{rankone}.


%









The preceding constructions depend on the parameters $\eta,\ep>0$, 
and we can make them
sufficiently small. Let $\ty_i,\,1\le i,$ be as above, and set
$\beta_i=\pi([\tx,\ty_i])$. By construction, $\beta_i\in\Ga(x,y)$,
and $|\beta_i|\to\infty$. Also by construction, $\beta_i$
belongs to $T_{\ep}(\alpha)$ except, possibly, for the interval of length
$\rho$ in the beginning and the interval of length
$\rho + r$ at the end. Hence, if the geodesic $\beta_i$ passes through
a point $f\in F$, it happens either during the first $\rho$ 
or during the last $\rho + r$ units of its lifespan.

For any pair of points $a,b\in M$, and any $l>0$ there is only a finite
number of geodesics with endpoints $a,b$ of length less than $l$.
Therefore, at most a finite number of geodesics $\beta_i,\,1\le i,$
can intersect $F$. The remaining infinite collection of connecting
geodesics $\beta_i$ does not pass through $F$. This
proves our claim in the rank one case.

\vspace{5mm}

From now until the end of this subsection, we assume that $\rk (M) \ge 2$.
Let $A \subset S$ be a $\Ga$-compact flat of maximal dimension; thus 
$\pi(A) \subset M$ is a maximal flat torus.
Unlike the rank one case, it is possible that  $F_1 = \pi(A) \cap (F\cup\{x,y\}) \ne \emptyset$.
Set $F_2 = (F\cup\{x,y\}) \setminus \pi(A)$, and let $\ep >0$ be such that
$T_{\ep}(\pi(A))  \cap F_2 = \emptyset$. Then $T_{\ep}(A)  \cap \tF_2 = \emptyset$;
let $\tx \in \pi^{-1}(x)\setminus A$. 

Let $g \in \Ga$ be a translation of $A$ in a regular direction. See
section~\ref{shad_sec} for details.
Let $\tp \in (A \setminus \tF_1)$ be an arbitrary point, and for $-\infty < i < \infty$
set $\tp_i = g^i\tp$, $\tz_i = s_{\tp_i}(\tx)$. There is a geodesic, $\alpha\subset A$, 
containing the points $\tp_i,\,-\infty < i < \infty$.
Let $0<\delta<\ep$ be such that $B_{\delta}(\tp) \cap \tF = \emptyset$. Then
$B_{\delta}(\tp_i) \cap \tF = \emptyset$ for all $i$. 

Let $0<\eta$. Set $\gamma_i=[\tx,\tz_i]$, and for arbitrary $\tw_i \in B_{\eta}(\tz_i)$ 
set $\gamma_{\tw_i}=[\tx,\tw_i]$. Note that $\gamma_i=\gamma_{\tz_i}$.
The regularity of $\gamma$ implies (see Lemma~\ref{hyperbolic_lem})
that just as in the rank one case, the finite geodesics $\gamma_{\tw_i}$
belong to $T_{\ep}(A)$ with the possible exception of their segments of
uniformly bounded length, located  in the beginning and at the end
of each geodesic.

By construction, for all $i$, the geodesic $\gamma_i$ intersects $A$ at a single point,
$\tp_i\notin\tF$. Apriori, a perturbed geodesic $\gamma_{\tw_i}$ could intersect
$A$ at a point of $\tF$, jeopardizing our proof. A geometric argument
based on considerations of symmetry and convexity, shows that if the intersection
$\gamma_{\tw_i}\cap A$ is nonempty, then it is stable under perturbations. 
See Lemma~\ref{nearhitting_lem}. By this lemma, $\{\gamma_{\tw_i}\cap A\}\subset B_{\delta}(\tp_i)$,
implying that the intersection point does not belong to $\tF$.
Note that the proof of Lemma~\ref{nearhitting_lem} crucially uses that $S$ is a 
symmetric space.

Extending the geodesics $\gamma_{\tw_i}$ beyond the point $\tw_i$
by $0<r<\infty$, we obtain the sequence of restricted shadows 
$SH_i(\tx,\eta,r)\subset S,\,1\le i$.
We claim that there exists $r$ such that for
all $i$ sufficiently large the projection $\pi:SH_i(\tx,\eta,r)\to M$ is surjective.
In a symmetric space of rank one the diameter of any restricted shadow $SH(\tx,\eta,r)$
grows with $r$ at a uniform rate. Due to the existence of multidimensional flats, this
fails in higher rank symmetric spaces. We will prove the claim using the uniform regularity
of the geodesics $\gamma_i,\,1\le i,$ and the  unique ergodicity
of the horocycle flow on compact locally symmetric spaces of noncompact type.
The latter is due to Hedlund \cite{He} for compact surfaces of constant
negative curvature,  and to Veech \cite{Ve} in the general case. See Theorem~\ref{HV_thm}
below. The
surjectivity follows from  Theorem~\ref{shadow_lem}, which
we call the {\em shadow lemma}. 
We deduce it from Corollary~\ref{cor_horocycle} of 
Theorem~\ref{HV_thm} and a uniform convergence. See Lemma~\ref{lem_convergence}.

The shadow lemma allows us to complete the proof of Theorem~\ref{nct_thm}
using the same argument as in the rank one case. Namely, we construct an infinite family
of connecting geodesics $\beta_i\in\Ga(x,y),\,i_0\le i$. Each geodesic
$\beta_i$ does not encounter points of $F$, except for, possibly, during the first $\rho$
or the last $\rho+r$ units of its life span. By preceding argument,
at most a finite number of the geodesics $\beta_i$ can pass  through $F$, 
contrary to the assumption that $F$ is a blocking set for $\{x,y\}$.



\subsection{Horocycles and the shadow lemma} \label{shad_sec}

Our proof of Theorem~\ref{nct} crucially uses a result about the density
of horocycles due to Hedlund and Veech. For convenience of
the exposition, we formulate it below.

\begin{thm} \label{HV_thm}
{\em (Hedlund, Veech)}

\noindent Let $M = \ {\Gamma} \setminus  S$ be a compact
locally symmetric space of noncompact type, and let $\pi:S\to M$ be the projection.

\noindent Let $\xi \in \partial_{\infty}S$ be regular, 
let $x \in S$ be arbitrary, and let $HC(\xi,x)\subset S$ be the corresponding
horocycle. Then $\pi(HC(\xi,x))$ is dense in $M$.
\end{thm}

\begin{rem} \label{He-Ve_rem}
{\em Hedlund \cite{He} proved Theorem~\ref{He-Ve_rem} for the 
hyperbolic plane. The general case follows from a theorem of
Veech \cite{Ve} about  the unique ergodicity of 
horocycle flows.}
\end{rem}

Let $HC(\xi,x)$ be a horocycle. For any $r>0$ we define the  
{\em restricted horocycle} $HC_r(\xi,x) = HC(\xi,x) \cap B_r(x)$.

\begin{cor} \label{cor_horocycle}
Let $M = \ {\Gamma} \setminus  S$  be as above.
For any $\ep > 0$ there exists
$r_0 = r_0(S,\Ga,\ep) > 0$ such that for
all $r > r_0$, any $x \in S$, and any regular point
$\xi \in \partial_{\infty}S$, the set
$\pi( HC_r(\xi,x))$ is $\ep$-dense in $M$.
\end{cor}

\begin{proof}
Let $x \in S$ and a regular point $\xi \in \partial_{\infty}S$ be given.
Let $HC(\xi,x)$ be the corresponding horocycle. Then
$\pi(HC(\xi,x))$ is dense, by Theorem~\ref{HV_thm}. Therefore,
there exist an open neighborhood $U\subset\partial_{\infty}S$ of $\xi$, 
an open neighborhood $V\subset S$ of $x$, and a positive number
$\rho = \rho(\xi,x,\ep,U,V)$, such that
$\pi(HC_{\rho}(\zeta,y))$ is $\ep$-dense in $M$ for
all $\zeta \in U$ and $y \in V$.
Using that $\Ga$ acts cocompactly, and Remark~\ref{rem_horocycle}, 
we obtain the claim.
\end{proof}

Let $\ga : \R\to S$ be a regular geodesic. Let
$\ga(\infty) = \xi \in \partial_{\infty}S$ and $\ga(0) = x$. 
Let $N \subset G$ be the nilpotent subgroup with
$HC(\xi,x) = N \cdot x$ as described in section
\ref{notcomp_set}. Note that then 
$HC(\xi,\ga(t)) = N \cdot \ga(t)$ for all
$t \in \R$.

\begin{lem} \label{lem_convergence}
Let $n \in N$.
Then the geodesics
$\ga(t)$ and $n \ga(t)$ converge
exponentially for $t\to \infty$.
The rate of convergence depends only on
$ \la_0 ^{+} (\ga)$.
\end{lem}

\begin{proof}
We will use two well known formulas from the
theory of Lie groups.
For $g \in G$ and $Y \in \g$ we have
$ g \grexp(Y) g^{-1} = \grexp(Ad(g)(Y))$;
for $X,Y \in \g$ we have $ Ad(\grexp(X))(Y) = e^{ad(X)}(Y) $.

We write
$\ga(t) = \grexp(tH) \cdot x$ and $n = \grexp(Y)$ with
$Y = \sum_{\la \in \Delta^{+}} Y_{\la}$. Then
\begin{align*} 
d(n\ga(t),\ga(t))  &= d(\grexp(Y)\grexp(tH) \cdot x, exp(tH) \cdot x) \\
&= d(\grexp(-tH)\grexp(Y)\grexp(tH) \cdot  x, x) \\
&= d(\grexp(Ad(\grexp(-tH))(Y)) x,x)  \\
&= d(\grexp(e^{ad(-tH)}(Y))\cdot x,x) \\
&= d(\grexp(\sum_{\la \in \Delta^{+}} e^{-t\la(H)}Y_{\la})\cdot x,x) 
\end{align*}

Since
$\la(H) \geq  \la_0^{+}(\ga) > 0$ for all
$\la \in \Delta^{+}$, the claim follows.
\end{proof}

\medskip

\begin{figure} \label{shadow}
\psfrag{xi}{$\xi$}
\psfrag{y}{$y$}
\psfrag{x}{$x$}
\psfrag{gamma}{$\gamma$}
\psfrag{gammar}{$\gamma(-r_{1})$}
\psfrag{ngamma}{$n \gamma$}
\psfrag{z}{$z$}
\psfrag{HCx}{$HC(\xi,x)$}
\psfrag{HCr}{$HC(\xi,\gamma(-r_{1}))$}
\includegraphics{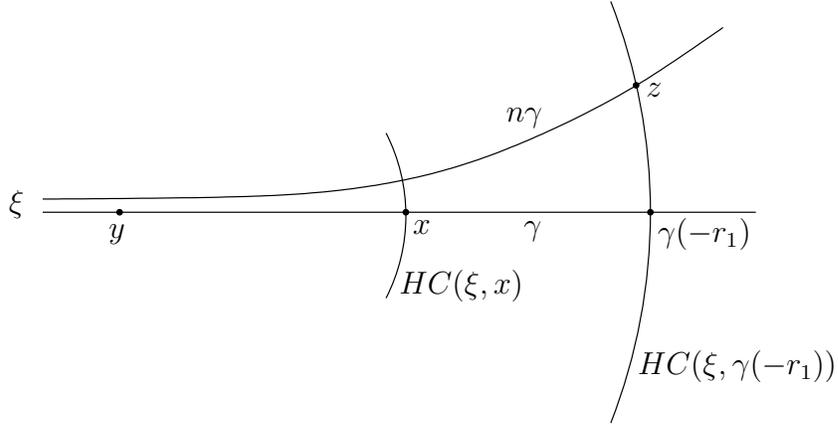}
\caption{Shadow Lemma}
\end{figure}

The following proposition, which is of independent interest, 
will be used in our proof of Theorem~\ref{nct}. For obvious reasons,
we call it the {\em shadow lemma}. See figure~\ref{shadow}.

\begin{thm} \label{shadow_lem}
Let $M=\Ga\setminus S$ be a compact, locally symmetric space of noncompact type. Then
for any $\ep,\eta > 0$ there exists
$R = R(S,\Gamma,\eta,\ep) >0$ so that the following holds.

\noindent Let $x,y \in S$ be distinct points. Suppose that the geodesic
$\ga$ containing them satisfies $\la_0^{+} (\ga) \geq \eta$.
Then the restricted shadow $SH(y,x,\ep) \cap B_R(x)$
has the property $\pi( SH(y,x,\ep) \cap B_R(x)) = M$.
\end{thm}

\medskip

\medskip


\begin{proof}
Let $l =d(x,y)$ and let
$\ga:\R \to S$ be the parametrization
with $\ga(0) =x$ and $\ga (l) =y$.
Let $\xi = \ga(\infty) \in \partial_{\infty}S$.
Consider as above the nilpotent subgroup
$N \subset G$ with
$HC(\xi,\ga(t)) = N\cdot \ga(t)$.
For $\sigma \ge 0$ and $t \in \R$ we define
$N(\sigma,t) \subset N$ to be the subset such that
$$HC_{\sigma}(\xi,\ga(t))= N(\sigma,t)\cdot \ga(t).$$
In other words
$N(\sigma,t)=\{n \in N: d(n\ga(t),\ga(t))\le \sigma\}$.

According to Corollary~\ref{cor_horocycle},
there exists
$r_0 = r_0(S,\Ga,\ep / 3)$ such that
$\pi(HC_{r_0}(\xi,\ga(t)))$ is
$\ep / 3$ dense in $M$ for all
$\xi \in \partial_{\infty}S$ and all $t\in \R$.

We claim that there exists
$r_1 = r_1 (\eta,r_0,\ep/3) >0$ such that
$$N(r_0,-r_1) \subset N(\ep/3,0).$$
To prove the claim, we assume that
$n \in N \setminus N(\ep/3,0)$.
This means that
$d(n\ga(0),\ga(0)) > \ep/3$.
By Lemma~\ref{lem_convergence}, there exists
$r_1$ depending only on $\eta, r_0$ and $\ep/3$, such that
$d(n\ga(-r_1),\ga(-r_1)) > r_0$.
Thus $n \in N\setminus N(r_0,-r_1)$. This proves our claim.

We set, for simplicity of notation,
$N_{\ep/3}= N(\ep/3,0)$.
The claim then implies
$$HC_{r_0}(\xi,\ga(-r_1)) \subset N_{\ep/3}\cdot \ga(-r_1).$$
Thus
for any $z \in HC_{r_0}(\xi,\ga(-r_1)$
there exists $n \in N_{\ep/3}$ such that
$z = n\cdot \ga(-r_1)$.
Then
$d(n\ga(0),x) \leq \ep/3$ and since
the function
$t \mapsto d(n\ga(t),\ga(t))$ is convex and thus
monotonously decreasing by Lemma~\ref{lem_convergence},
we have
$d(n \ga(l),y) \leq \ep/3$.
By the triangle inequality and the convexity of the
distance function,
$[y,z] \cap B_{2\ep/3}(x) \neq \emptyset$.
Thus
$$ HC_{r_0}(\xi,\ga(-r_1)) \subset SH(y,x,2\ep/3) .$$
By convexity again
$$ B_{\ep/3}(HC_{r_0}(\xi,\ga(-r_1))) \subset SH(y,x,\ep). $$
The triangle inequality also implies that
$$ B_{\ep/3}(HC_{r_0}(\xi,\ga(-r_1))) \subset B_R(x),$$
where $R = (r_1+ r_0 + \ep/3)= R(S,\Gamma,\eta,\ep)$. 
Hence
$$ B_{\ep/3}(HC_{r_0}(\xi,\ga(-r_1))) \subset SH(y,x,\ep) \cap B_R(x).$$ 
Since $\pi(HC_{r_0}(\xi,\ga(-r_1)))$ is $(\ep/3)$-dense in $M$, 
this implies 
$$\pi(SH(y,x,\ep)\cap B_R(x)) = M.$$
\end{proof}



\medskip

\subsection{Proof of Insecurity} \label{proof_insec}
In this section we prove Theorem~\ref{nct_thm}.
Let $M = \  {\Gamma} \setminus S $ be a compact,
locally symmetric space of noncompact type, and
let $\pi: S \to M$ be the covering map.
Let $A \subset S$ be a
$\Ga$-compact flat. (They
are dense in the set of all flats \cite[Lemma 8.3]{Mo}.)
Then $A$ is 
totally geodesic and isometric to
$\R^{\rk (S)}$.
Since $A$ is $\Gamma$-compact,  there exists a subgroup
$\Ga_A \sim \Z^{\rk (S)}$ of $\Ga$
which operates by translations with a compact quotient on $A$.
We also choose a point
$\tx \in \pi^{-1}(x)$ with $\tx \not\in A$
and a point
$\tp \in A$ such that $\tp \not\in \tF$,
where
$\tF = \pi^{-1}(F)$.
Since $\pi(A)$ is a closed subset of $M$, there
exists a constant
$\ep_1 > 0$ such that $(T_{2\ep_1}(A) \setminus A) \bigcap \tF = \emptyset$.
Without loss of generality, we assume that
$4 \ep_1 < d(\tx,A)$.
Let $g \in \Ga_A$ be a translation in a regular
direction. Hence, for any $q \in S$ the point
$\lim_{i\to\infty} g^i(q) \in \partial_{\infty}(S)$ is regular.
For $i \in \N$ set $\tp_i = g^i\tp \in A$.
Since $A$ and $\tF$ are $g$-invariant, and
$\tp \not \in \tF$, there exists
$\ep_2 >0$ such that
$B_{\ep_2}(\tp_i) \cap \tF = \emptyset$.
We assume without loss of generality that $\ep_2 \leq \ep_1$.
Set $\tz_i = s_{\tp_i}(\tx)$, and  $b_i = d(\tx,\tp_i)$. Let
$\ga_i:\R \to S$ be the geodesic determined by
$\ga_i(0) = \tx$ and $\ga_i(b_i) = \tp_i$. Then $\ga_i(2b_i)=\tz_i$.

\vspace{3mm}

We will now state two lemmas, and derive Theorem~\ref{nct_thm} from them.
Then we will prove the lemmas.

\begin{lem} \label{hyperbolic_lem}
There exists $\rho > 0$ such that
for all $i \in \N$ we have
$\ga_i([\rho,b_i]) \subset T_{\ep_1}(A)$.
\end{lem}

\begin{lem} \label{nearhitting_lem}
There exists $\ep_3 > 0$ such that for every
point $\tw_i \in B_{\ep_3}(\tz_i)$
either
$[\tx,\tw_i] \cap A = \emptyset$ or
$[\tx,\tw_i] \cap A$ is a point whose distance from $\tp_i$ is less
than $\ep_2$. 
In particular, $[\tx,\tw_i] \cap A \cap \tF = \emptyset$.
\end{lem}

\noindent {\em Proof of Theorem~\ref{nct_thm}.}
Let $\ga_i$ be the geodesics  defined above with
$\ga_i(0) = \tx$,
$\ga_i(b_i) = \tp_i$ and
$\ga_i(2b_i)=\tz_i$.
By Lemma~~\ref{hyperbolic_lem},
$\ga_i([\rho, b_i ]) \subset T_{\ep_1}(A)$
and 
using the geodesic symmetry at $\tp_i$, we see
that
$\ga_i([\rho,2 b_i -\rho]) \subset T_{\ep_1}(A)$.
Note that the geodesics
$\ga_i$ converge to a limit geodesic
$\ga_{\infty}$, with
$\ga_{\infty}(0) = \tx$ and
$\ga_{\infty}(\infty)= \lim_{i\to\infty} g^i(\tp) \in \partial_{\infty}(S)$.
Since $\ga_{\infty}(\infty)$ is regular,
$ \la_0^{+} (\ga_{\infty}) > 0$.
Thus, by passing to a subsequence, if necessary, we insure
that there exists $\eta > 0$ such that
$ \la_0^{+} (\ga_i) \geq \eta$ for all $i$.
We will use the number $\ep_3$ of Lemma \ref{nearhitting_lem}, assuming,
without loss of generality, that $\ep_3 \leq \ep_1$. Let
$R = R(S,\Ga,\eta,\ep_3/2) > 0$ be the
number from Theorem~\ref{shadow_lem}.

Then there exists
$\ty_i \in B_R(\tz_i) \bigcap SH(\tx,\tz_i,\ep_3 /2)$ 
with $\pi(\ty_i) = y$.
Let
$\sigma_i:[0,\ell_i] \to \tM$ be the unit speed
parametrization of the geodesic
$[\tx,\ty_i]$. Then $\ell_i \leq 2 b_i +R + \ep_3$.
Since $[\tx,\ty_i]$  comes $\ep_3 /2$ close
to $\tz_i$, the triangle inequality implies
$d(\sigma_i(2 b_i),\ga_i(2 b_i)) \leq \ep_3$.
By convexity, $d(\sigma_i(t),\ga_i(t)) \leq \ep_3$
for $t \in [0,2 b_i]$, and since $\ep_3 \leq \ep_1$, we have
$\sigma_i([\rho,2 b_i -\rho]) \subset T_{2 \ep_1}(A)$.
By Lemma~\ref{nearhitting_lem}, 
$\sigma_i$ does not intersect $A$ at a point of $\tF$,
thus $\sigma_i([\rho, 2 b_i - \rho]) \cap \tF = \emptyset$.
Hence, for $\rho' = \rho + R + \ep_3$ we have
$\sigma_i([\rho,\ell_i - \rho']) \cap \tF = \emptyset$.

The geodesics  $\tau_i = \pi \circ \sigma_i :[0,\ell_i] \to M$
connect $x$ and $y$. By a discreteness argument, there are at most 
finitely many indices $i$ such that
$\tau_i([0,\rho])  \cap F \neq \emptyset$ or
$\tau_i([\ell_i-\rho',\ell_i])\cap F \neq \emptyset$.
We have thus constructed infinitely many connecting geodesics that do not meet $F$.
\qed

\medskip

\noindent {\em Proof of Lemma~\ref{hyperbolic_lem}.} 
Let $\omega = d(\tp,g\tp)$, and let $\ga: \R \to A$ be the unit speed geodesic
with $\ga(0) = \tp$ and $\ga(\omega) = g\tp$.
Then $\ga$ is the axis of the isometry
$g$ passing through the point $\tp$. Since
$g$ is a translation in a regular direction,
$\ga$ is a regular geodesic, i.e.,  all
parallels to $\ga$ are contained in  $A$.
For $i\in\N$ we have $\ga(i\omega) = g^i\tp$.
Set $c = d(\tx,\tp)$.
Then $d(\ga_i(0),\ga(0))=c$, and by
triangle inequality,
$d(\ga_i(b_i),\ga(b_i))\leq c$; thus by convexity,
$d(\ga_i(t),\ga(t))\leq c$ for all $t \in [0,b_i]$.
For $i\in\N$ set $r_i:[0,b_i] \to [0,\infty)$ be the 
function $r_i(t) = d(\ga_i(t),A)$.
Since the curvature of $S$ is nonpositive, these functions
are convex, and since $r_i(b_i) = 0$  they are nonincreasing.
In view of $r_i(0) = d(\tx,A) \geq 3 \ep_1$, 
there exists a unique $\rho_i \in (0,b_i)$ such that
$r_i(\rho_i) = \ep_1$.

It remains to show that there exists a $\rho$ such that
for all $i \in \N$ we have $\rho_i \leq \rho$. 
Assume the opposite. Then, passing to a subsequence, if necessary,
we have $\rho_i \to \infty$.
Let $s_i = \rho_i / 2$ and let
$j(i) \in \N$ be such that
$ \mid j(i) \omega - s_i \mid \leq \frac{\omega}{2}$.
Then $d(\ga_i(s_i),\tp_{j(i)})\leq c+\frac{\omega}{2}$.
We reparametrize the geodesics
$g^{-j(i)}\circ \ga _{\mid [0,\rho_i]}$ so that
the parameter interval is $[-s_i,s_i]$, thus obtaining
$\ga^*_i:[-s_i,s_i] \to A$, where
$\ga^*_i(t) = g^{-j(i)}\circ \ga (t + s_i)$.
By construction,
$d(\ga^*_i(t),\ga(t)) \leq c+\frac{\omega}{2}$ and 
$\ep_1\leq d(\ga^*_i(t),A)\leq d(\tx,A)$.
Hence, there is a converging subsequence
$\ga^*_i \to \ga^*$, the limit geodesic $\ga^*$ is defined on $\R$. 
The function $d(\ga^*(t),\ga(t))$ on $\R$ is convex and bounded,
hence $d(\ga^*(t),\ga(t))=d$ is constant, i.e., $\ga^*$ is parallel to $\ga$;
since $d\neq 0$, the geodesic $\ga^*$ is not contained in $A$. This  contradicts to the
regularity of $\ga$.\qed

\medskip

\noindent {\em Proof of Lemma~\ref{nearhitting_lem}.} 
We will show that
$\ep_3 = \min \{ \frac{\ep_1 \ep_2}{2 \rho},\frac{\ep_2}{4} \}$
satisfies the requirements.
Let $\tw_i \in B_{\ep_3}(\tz_i)$ and assume that
$[\tx,\tw_i]\cap A \neq \emptyset$.
Let
$\ga_i: \R \to \tM$ be the unit speed geodesic
with
$\ga_i(0) = \tx$ ,$\ga_i(b_i) = \tp_i$ and
$\ga_i(2b_i) = \tz_i$.
Let $\sigma_i:\R \to \tM$ be the unit speed geodesic
with $\sigma_i(0)= \tx$ and $\tw_i = \sigma_i(d_i)$
for some $d_i > 0$.
Since $\ga_i(2 b_i) = \tz_i$, by construction,
$d(\sigma_i(2 b_i),\ga_i(2 b_i)) \leq 2 \ep_3 $, and 
by convexity,
$d(\sigma_i(t),\ga_i(t)) \leq 2 \ep_3$  for $t \in [0,2 b_i]$.
Since $[\tx,\tw_i]\cap A \neq \emptyset$,
there exists a unique point $h_i \in [0,2 b_i]$ 
 satisfying $\sigma_i(h_i) \in A$.
It suffices to show that
$\mid h_i - b_i \mid \leq \frac{\ep_2}{2}$.
Indeed, then
$d(\sigma_i(h_i),\tp_i) =
d(\sigma_i(h_i),\ga_i(b_i))
\leq  \frac{\ep_2}{2}  + 2\ep_3 \leq \ep_2$.

To prove the inequality $\mid h_i - b_i \mid \leq \frac{\ep_2}{2}$,
we will use the functions
$f_i(t) = d(\sigma_i(t),A)$ on $\R$.
Since $A$ is invariant under the geodesic
reflection at the point $\sigma_i(h_i)$,
the function
$f_i$ is symmetric with respect to
$h_i$.

The functions $f_i$ have the following properties:

\begin{enumerate}
\item  $f_i$ is nonnegative, convex, $1$-lipschitz, and symmetric
with respect to $h_i$;
\item $f_i^{-1}(0) =\{ h_i\}$ and for $t > 0$ the set $f_i^{-1}(t)$ consists of two elements;
\item $f_i(0) \geq 4\ep_1$;
\item $f_i(\rho) \leq 2\ep_1$;
\item $\mid f_i(0) -f_i(2b_i) \mid \leq 2 \ep_3$.
\end{enumerate}

Items (1)-(3) are obvious. To show (4), we note that
$d(\ga_i(\rho), A) \leq \ep_1$ and
$d(\sigma_i(\rho),\ga_i(\rho))\leq 2 \ep_3\leq \ep_1$;
(5) follows from $f_i(0) =d(\tx,A)=d(\tz_i,A)=d(\ga_i(2b_i),A)$, and
$d(\sigma_i(2b_i),\ga_i(2b_i))\leq 2 \ep_3$.

By (1), (3) and (4), there is
$c_i \in [\ep_1,\rho]$ such that
$f_i(c_i) = 3 \ep_1$.
By convexity, we have $ f_i'(t) \leq - \frac{\ep_1}{\rho}$
for $0 \leq t\leq c_i$. Hence, for $0 \leq t \leq \ep_1$, we have the
inclusion $f_i([-t,t]) \supset [f_i(0) - t \frac{\ep_1}{\rho},f_i(0) + t \frac{\ep_1}{\rho}]$.
By (5), $f_i(2 b_i) \in f_i([-t,t])$ for
$t \geq \frac{2 \ep_3 \rho}{\ep_1}$.
Thus, there exists $t_i,\,\mid t_i \mid \leq \frac{2 \ep_3 \rho}{\ep_1}$,
such that $f_i(t_i) = f_i(2 b_i)$.
By symmetry, $h_i = \frac{1}{2}(2 b_i + t_i)$, implying
$\mid h_i - b_i \mid = \frac{1}{2}\mid t_i \mid \leq \frac{ \ep_3 \rho}{\ep_1} 
\leq \frac{\ep_2}{2}.$ \qed


\newpage


\section{Locally symmetric spaces of compact type}
\label{comp_type}
We first recall the basics on symmetric spaces of  compact  type.
\subsection{Symmetric spaces of compact type}       
\label{comp_set}

A (simply connected) symmetric space of compact type satisfies $S=G/K$, where
$G = \iso_0(S)$ is a compact, connected, semisimple Lie group, 
and $K\subset G$ is (essentially) the fixed point set of an involution
$\sigma:G\to G$. The material of section~\ref{basic_sub} and much of
that of section~\ref{notcomp_set} applies, and we will use the
notation established there. In contrast to section~\ref{notcomp_set},
it is convenient to fix once and for all a reference point, $o\in S$. 
Hence, we also fix the Cartan decomposition $\g = \k + \p$,
and the identification $\p = T_oS$.
If $p:G\to S$ is the projection, then $o=p(e)$ and $K=p^{-1}(o)$. Let
$s:S\to S$ be the geodesic involution about $o$. Then $\sigma(g)=sgs$.

Our exposition does not depend on a $G$-invariant riemannian metric 
on $S$; for concreteness, we choose the metric $<\cdot,\cdot>$ corresponding to the
negative of Killing form. As in section~\ref{notcomp_set}, $\riexp:\p\to S$ and
$\grexp:\p\to G$ are the riemannian and the Lie group exponential maps
respectively. They satisfy  
$\riexp X = p(\grexp X),\ \riexp(\adj_{\p}k(X)) = p(k\cdot\grexp X)$.

A {\em (maximal) flat}, $B\subset S$, is a (maximal)
totally geodesic submanifold, isometric to a flat torus.
The {\em rank} $\rk(S)$ is the   dimension of a maximal flat. 
Let $\A$ be the set of {\em maximal abelian subalgebras} in $\p$,
and let $\ttt_o$ be the set of maximal flats in $S$ containing $o$.
The mapping $\riexp$ yields a $K$-equivariant isomorphism between $\A$ and $\ttt_o$;
therefore the action of $K$ on $\ttt_o$ is transitive.
Hence, all  maximal flats in $S$ are isometric to a 
flat torus of dimension $\rk(S)$.
We will refer to the flats $B\in \ttt_o$ as the {\em maximal tori} in $S$. 
We fix a {\em reference Cartan subspace}, $\a\subset\p$, and let
$A=\riexp(\a)$ be the {\em reference maximal torus}. (We also denote by $A\subset G$ 
the corresponding subgroup.) 

For $x\in S$ (resp. $X\in\p$) let $\ttt_o(x)\subset\ttt_o$ (resp. $\A(X)\subset\A$) 
be the set of maximal tori (resp. Cartan subspaces) containing $x$ (resp. $X$).
A point $x\in S$ (resp. a vector $X\in\p$) is {\em regular} if its $K$-isotropy
subgroup  $K_x\subset K$ (resp. $K_X\subset K$) has minimal dimension. 
Then $x\in S$ (resp. $X\in\p$) is regular iff 
the set $\ttt_o(x)$ (resp. $\A(X)$) consists of a single element.
The set $S_r$ of regular points in $S$ (resp. the set $\p_r$
of regular  vectors in $\p$)
is $K$-invariant. We refer to the elements of the complement
$S\setminus S_r$ (resp. of the set $\p\setminus\p_r$) as {\em singular}
points. The {\em singular set} $S\setminus S_r$ has a natural stratification. 
The maximally singular points $x\in S\setminus S_r$ satisfy $K_x=K$. 
The finite set $Z=Z(S)$ of these points is the
{\em center} of the symmetric space \cite{Lo}. It depends on the
reference point. Note that $o\in Z$, and $Z=\cap_{B\in\ttt_o}B$ \cite{Lo}.

\medskip

Let $X$ be a riemannian manifold, and let $Y\subset X$ be a closed submanifold. If
$x,y\in Y$ we denote by $\Gamma_Y(x,y)\subset\Gamma(x,y)$ the subcollection of 
connecting geodesics that belong to $Y$.
\begin{lem}  \label{algebr_lem}
Let $x,y\in S$ be arbitrary points, and let $\gamma\in\Gamma(x,y)$.
Then there exists a maximal torus $B \subset S$ such that $\gamma\subset\Gamma_Y(x,y)$.
\end{lem}  
{\bf Proof}. By homogeneity, it suffices to establish the claim for
$\gamma\in\Gamma(o,x)$, where $x\in S$ is an arbitrary point.
By preceding remarks, $\gamma=\{\riexp\, tX, 0\le t\le |\gamma|\}$, 
for some $X\in\p$.
The vector $X$ is contained in a Cartan subspace $\b\subset\p$. 
But the maximal torus $B=\riexp(\b)$ is totally geodesic. \qed

\subsection{Group action and security}   \label{action}
We will investigate  the security
of manifolds with a group action. In section~\ref{block_tori} we will
apply this material to symmetric spaces of compact type.

Let $M$ be a compact riemannian manifold. (It is not, in general, a symmetric
space). Let $U\subset\iso(M)$ be a closed, infinite subgroup.
For $x\in M$ we denote by $U\cdot x\subset M$ the $U$-orbit of $x$, and
by $U_x \subset U$ the isotropy subgroup of $x$. 

Let $\gamma$ be a geodesic in $M$, and let $z\in\gamma$.
If $T_z(U\cdot z)\subset T_zM$ is orthogonal to 
$\gamma$ at $z$ we will say that $U$ {\em acts transversally to $\gamma$
at the point $z$}. The following is motivated by Definition 2.1 in \cite{BS}.  

\begin{defin}  \label{trans_def}
A geodesic $\gamma\subset M$ is {\em transversal} (to the action
of $U$) if  $U$ acts transversally to $\gamma$
at any point $z\in\gamma$. A {\em collection $\Gamma$ of geodesics in $M$
is transversal} (to the action of $U$) if every $\gamma\in\Gamma$
is transversal.
\end{defin} 

A group acting on a manifold, naturally acts on the set of configurations.
Extending the notation above, we
will denote this action by $u\cdot\{x,y\}$. Then $U\cdot\{x,y\}$ is the $U$-orbit
in the space of configurations.
A configuration $\{x,y\}$ is {\em fixed} if $U\cdot\{x,y\}=\{x,y\}$.
If $U$ is connected, then $\{x,y\}$ is fixed iff both $x,y$ are fixed points.

\begin{prop}  \label{act_block_prop}
Let $M$ and $U$ be as above; let $\{x,y\}$ be a configuration such that 
$\Gamma(x,y)$ is a transversal collection of geodesics.\\
If $U^0$ fixes the  configuration $\{x,y\}$,
then $\{x,y\}$ is secure iff it has a blocking set that
consists of $U^0$-fixed points.
\end{prop}
{\bf Proof}. We will assume, for convenience of exposition, that $U$ is connected
and that the isotropy subgroups $U_z$ are connected for all $z\in M$. The general 
situation reduces to this case case by passing to the identity components
of relevant groups. We normalize the double-invariant riemannian metric
of $U$ such that the mappings $g\mapsto g\cdot z$, etc do not increase the relevant distances.
Let $X$ be a complete, riemannian manifold, and let $Y\subset X$
be an arbitrary subset. For $r>0$ and $z\in X$
we denote by $Y(r,z)$ the intersection of $Y$ with the open ball of radius 
$r$ in $X$ centered at $z$.

Let $N\subset M$ be the set of $U$-fixed points.
Then $x,y\in N$ and the collection $\Gamma(x,y)$ is $U$-invariant.
Let $B$ be a minimal blocking set. It suffices to show that $B\subset N$.
Suppose that this fails, and set $B_0=B\cap N,\,B_1=B\setminus B_0$.
By assumption, $B_1\ne\emptyset$.
By minimality of $B$, there exists $\gamma\in\Gamma(x,y)$
such that $\gamma$ does not pass through $B_0$. 
Let $z_1,\dots,z_m$ be the points of $B_1$ contained in $\gamma$,
and let $U_1,\dots,U_m\subset U$ be their isotropy subgroups.
They are proper subgroups of $U$, hence 
$\Omega=U\setminus(U_1\cup\dots\cup U_m)\subset U$ is a dense open set.

We define the mapping $\varphi:U\to M^m$ (the $m$-fold product)
by $\varphi(g)=(g\cdot z_1,\dots,g\cdot z_m)$. Denote by $X\subset M^m$
the subset given by conditions
$$
X=\{(w_1,\dots,w_m):w_i\ne z_i,\,1\le i \le m.\}
$$
Then $\varphi$ is a differentiable map, 
$\varphi(U(e,\ep))\subset M^m((z_1,\dots,z_m),\ep)$, and
$\varphi(\Omega\cap U(e,\ep))\subset M^m((z_1,\dots,z_m),\ep)\cap X$.

Let $0\le t \le |\gamma|$ be the natural parameter, and let
$0<t_1<\cdots<t_m<|\gamma|$ be given by $\gamma(t_i)=z_i,\,1\le i \le m.$
Let $\ep>0$ be arbitrary.
For any $g\in U(e,\ep)$ the geodesic $g\cdot\gamma\in\Ga(x,y)$
is $\ep$-close to $\gamma$ (pointwise). Since $B$ is a finite set,
the distance $\delta=d(\gamma,B\setminus\{z_1,\dots,z_m\})>0$.
Hence if $\ep<\delta/2$ and $g\in \Omega\cap U(e,\ep)$, then 
the geodesic $g\cdot\gamma$ does not pass
through the points of $B\setminus\{z_1,\dots,z_m\}$.
By preceding remarks, if $0<\ep$ is sufficiently small
and $g\in \Omega\cap U(e,\ep)$, then $g\cdot\gamma$ does not pass 
through the points $z_i,\,1\le i \le m,$ either.
Thus, for any $g\in \Omega$ and sufficiently close to the identity,
the geodesic $g\cdot\gamma\in\Ga(x,y)$  is not blocked by $B$.
Hence, contrary to the assumption, $B$ is not a blocking set. \qed

\medskip

The following consequence of Proposition~\ref{act_block_prop} will be useful.

\begin{corol}  \label{act_block_cor}
Let $M$ be a compact riemannian manifold, and let $U$ be a compact Lie 
group of isometries of $M$.
Denote by $N\subset M$ the set of $U^0$-fixed points.

Let $\{x,y\}$ be a configuration  such that $\Gamma(x,y)$ is transversal 
to the action of $U$.\\
1. Let $x,y\in N$. Suppose that $\{x,y\}$ is a secure configuration, and let
$B\subset M$ be a blocking set. Then $B\cap N$ is a blocking set as well.\\
2. Let $|U\cdot\{x,y\}|<\infty$.
Then the configuration $\{x,y\}$ is secure iff it has a $U$-invariant
blocking set.
\end{corol}
{\bf Proof}. The first claim was actually obtained in the proof of
Proposition~\ref{act_block_prop}.
Set $X=(U\cdot x)\cup(U\cdot y)$. The set $X$ is finite, and, by
\cite{Gut04}, $\{x,y\}$ is secure iff the collection $\Gamma(X)$ of geodesics
connecting the points of $X$ has a finite blocking set.
Since $X\subset N$, the second claim now follows from the first. \qed

\medskip

Let $M$ be as above, and let $K$ be an infinite, compact group,
{\em properly acting} on $M$ by isometries. (Equivalently,
$K\subset\iso(M)$.) We will say that a geodesic $\gamma$
(resp. a collection $\Gamma$ of geodesics) in $M$ is {\em $K$-transversal}
if $\gamma$ (resp. any $\gamma\in\Gamma$) satisfies the requirements
of Definition~\ref{trans_def}.

The following theorem connects the preceding material with our main subject.
\begin{thm}  
\label{act_thm}
Let $M$ be a compact riemannian manifold, and let $K$ be an 
infinite, compact group, properly acting on $M$ by isometries.
Let $F\subset M$ be the set of $K_0$-fixed points.\\
Suppose that $F$ is a finite, nonempty set. 
Then the manifold $M$ is not secure. 
\end{thm}
{\bf Proof}. Let $\gamma\subset M$ be any geodesic.
By, Proposition~2.2 of \cite{BS}, $\gamma$ is $K$-transversal iff
$K$ acts transversally to $\gamma$ in at least a point, $z\in\gamma$.
(The point in question may be an endpoint of $\gamma$  as well.)
The transversality condition is trivially satisfied if $z\in F$.
Hence, any geodesic intersecting $F$ is $K$-transversal.

Let $\{x,y\}$ be a configuration such that $\{x,y\}\cap F\ne\emptyset$.
By preceding remarks, $\Gamma(x,y)$ is $K$-transversal.
We will now specialize to configurations $\{x,y\}$ such that both $x,y\in F$,
and consider two cases.

\noindent 1. Let $|F|=1$, and set $F=\{x\}$. By Proposition~\ref{act_block_prop},
the configuration $\{x,x\}$ is insecure.

\noindent 2. Let $|F|>1$, and let $x\in F$ be arbitrary.
Let $y\in F$ be a point, different from $x$, and such that the distance
$d(x,y)$ is less than or equal to $d(x,y'),y'\in F$,
for all $y'\ne x$. By Proposition~\ref{act_block_prop},
the configuration $\{x,y\}$ is secure iff every geodesic in $\Gamma(x,y)$
passes through $F$.
Let $\gamma\in\Gamma(x,y)$ be a geodesic such that $|\gamma|=d(x,y)$.
By construction, $\gamma$ does not pass through $F$.

Thus, in both cases $M$ has an insecure configuration. \qed

\subsection{Secure and insecure configurations}  \label{block_tori}
We will use the notation of section~\ref{comp_set}. 
Let $S=G/K$ be a symmetric space of compact type, and let $o\in S$ be the
reference point. Following \cite{BS}, we classify the points of $S$ by dimensions
of their $K$-orbits. Regular points  $x\in S$ are such that
$K\cdot x$ has the maximal dimension, $r=r(S)$.
The {\em defect of a point} is defined by $\delta(x)=r(S)-\dim(K\cdot x)$. 
See Definition~7.1 in
\cite{BS}. Thus, $x\in S$ is singular iff $\delta(x)>0$. The maximal
possible defect is $r(S)$, and the set of maximally singular points is
the center $Z\subset S$. 
The $K$-stabilizer, $L\subset K$, of a regular point is determined up
to conjugacy, hence $l(S)=\dim(L)$ is well defined. Note that
$\dim(K)=l(S)+r(S)$. Let now $\{x,y\}\in C(S)$ be a configuration.
The mapping $\{x,y\}\mapsto G_x\cap G_y$ is equivariant with respect
to the natural actions of $G$. Let $g \in G$ satisfy $g\cdot x = o$, 
and set $g\cdot y = w$. Then
$$
\dim(G_x\cap G_y)=\dim(K\cap G_w)= \dim(K_w).
$$
\begin{defin}  \label{reg_def}
The {\em defect of a configuration} is given by
\begin{equation}    \label{def_eq}
\delta(\{x,y\})=\dim(G_x\cap G_y)-l(S).
\end{equation}
A configuration $\{x,y\}$ is {\em regular}
if $\delta(\{x,y\})=0$, and {\em singular} if $\delta(\{x,y\})>0$.
\end{defin} 

We formulate the basic properties pertaining  to the regularity of configurations
in the proposition below.
\begin{prop}    \label{reg_prop}
Let $S=G/K$ be a symmetric space of compact type,
let $o\in S$ be the reference point, and let $A$ be a reference
torus. Let $\{x,y\}\in C(S)$ be arbitrary.

\noindent 1. The defect of a configuration is invariant with respect 
to the action of $G$ on $C(S)$. We have
$$
0 \le \delta(\{x,y\}) \le r(S).
$$

\noindent 2. The configuration $\{x,y\}$ is regular iff there is a unique maximal torus
containing $x,y$ iff $\{x,y\}$ is conjugate to $\{o,a\}$ where
$a\in A$ is a regular point.

\noindent 3. The configuration $\{x,y\}$ is singular iff
$\{x,y\}$ is conjugate to $\{o,a\}$ where
$a\in  A$ is a singular point.

\noindent 4. The equality $\delta(\{x,y\}) = r(S)$ holds iff 
$G_x=G_y$ iff  $\{x,y\}$ is conjugate to $\{o,z\}$ where
$z\in Z$.
\end{prop}
{\bf Proof}. The claims readily follow from the preceding discussion and 
the standard material \cite{He,He1,Lo}. \qed

\medskip

If $\delta(\{x,y\}) = r(S)$, we will say that the 
configuration $\{x,y\}$ is {\em maximally singular}.

\begin{thm}    \label{block_thm}
Let $S$ be a symmetric space of compact type. 

\noindent 1. Any regular configuration in $S$ is secure;
it has a blocking set of $2^{\rk(S)}$ points.

\noindent 2. There exist maximally singular configurations in $S$
that are insecure.
\end{thm}
{\bf Proof}. 1. Let $\{x,y\}$ be a regular configuration,
and let $B$ be the unique maximal torus containing $x,y$. 
See Proposition~\ref{reg_prop}. By Lemma~\ref{algebr_lem},
$\Ga(x,y)=\Ga_B(x,y)$. By Proposition~\ref{euclid_case}, a flat torus of $r$ dimensions
is uniformly secure; its security threshold is $2^r$.

\noindent 2. In view of  Proposition~\ref{reg_prop}, it suffices to consider 
the configurations $\{x,y\}$, where $x,y\in Z$. 
Our setting satisfies the assumptions of Theorem~\ref{act_thm} (with $F=Z$),
which implies the claim. \qed

\medskip

\begin{corol}    \label{block_cor}
Let $M$ be a locally symmetric space of compact type. 
Then almost all configurations in $M$ are secure.
However, $M$ always has insecure configurations as well.
\end{corol}
{\bf Proof}. The space $M$ has a finite covering $q:S\to M$, where
$S$ is a symmetric space of compact type, and $q$ is a local isometry.
The Lebesgue measure on the space of configurations in $M$ is the image under
$q_*$ of the corresponding  measure in $S$. Regular configurations in $S$
form a subset of full measure. The claims now follow from Theorem~\ref{block_thm},
and the basic facts concerning  security and coverings \cite{Gut03}. \qed
 
\medskip

We will now illustrate preceding propositions with a few examples.
First, we consider the case of symmetric spaces of rank one.

\subsection{Example: Security for symmetric spaces of rank one}  \label{exam1_sub}
Let $S=G/K$ be a symmetric space of rank one, and let $o\in A\subset S$ 
be the reference point and the reference torus.
Then $A$ is a circle; let $a'$ be the antipodal point of $a\in A$. 
A point $x\in S$ is either regular or maximally singular.
There are two possibilities: $Z=\{o\}$, or $Z=\{o,o'\}$.
Theorem~\ref{block_thm} and Proposition~\ref{reg_prop} show
that in the former case the configuration $\{o,o\}$ is insecure,
while in the latter case $\{o,o\}$ is secure, and $\{o,o'\}$ 
is insecure.

Let us investigate the security of general configurations in $S$. If $|Z|=2$
then the involution $a\mapsto a'$ extends by homogeneity to all of $S$.
Thus, we have the antipodal involution $x\mapsto x'$, and it satisfies
$g\cdot x' = (g\cdot x)'$. Then the {\em antipodal configurations} 
$\{x,x'\}$ are insecure; all other configurations are secure.
If $|Z|=1$ then the antipodal involution on $S$ does not exist.
The insecure configurations are $\{x,x\}$, and all other configurations are secure.

Compact symmetric spaces of rank one are listed in \cite{He2}, p. 518 and p. 535.
We will now illustrate the preceding discussion by briefly  going   
over the list.

\medskip

\noindent{\em The round sphere $S^n,\,n>1$}. 
(The space $S^1$ is of euclidean type.) We have
$|Z|=2$. The notion of antipodal points $x',x\in S^n$ is classical.
The only insecure configurations are 
$\{x,x'\}$ where $x$ is arbitrary. All of the geodesics in the continuum $\Ga(x,x)$
are blocked by $x'$.

\medskip

\noindent{\em The  real projective space $\R P^n,\,n>1$}.
(Note that $\R P^1=S^1$.)
Here $|Z|=1$, and the insecure configurations are
$\{x,x\}$ where $x\in S^n$ is arbitrary. For all  configurations
$\{x,y\}$ with $x\ne y$ the set $\Ga(x,y)$ is finite. Note that 
$S^n$ is a double covering of $\R P^n$. Denote the antipodal involution by 
$\sigma:S^n\to S^n$. Then $\R P^n=S^n/\sigma$. The remarks above illustrate 
Corollary~\ref{cov_sec_cor}.

\medskip

\noindent{\em The  complex projective space $\C P^n,\,n>1$}.
(The case $n=1$ is exceptional, and $\C P^1=S^2$.)
Now $|Z|=1$, hence the only insecure configurations are
$\{x,x\}$ where $x\in \C P^n$ is arbitrary.
Just like for $\R P^n$, only for these configurations the 
set of connecting geodesics is infinite. Note that $\C P^n$
is simply connected, and $\dim(\C P^n)=2n$.

\medskip

\noindent{\em The  quaternionic projective space $\HH P^n,\,n>1$}.
(The case $n=1$ is already disposed of, since $\HH P^1=S^4$.)
The space is simply connected, $\dim(\HH P^n)=4n$.
Again, $|Z|=1$, and the only insecure configurations are
$\{x,x\}$ where $x\in \HH P^n$ is arbitrary.
Just like for $\R P^n, \C P^n$, these are the only configurations for which the 
set of connecting geodesics is infinite.

\medskip

\noindent{\em The unique exceptional compact symmetric space
of rank one: The space $F\ II$}. See \cite{He2} (pp. 516, 518) and \cite{W}.
The isometry group of the space $F\ II$ is $F_4$: The compact,
connected Lie group with the Lie algebra $\f_4$. We have $\dim(F\ II)=16$
and $\dim(F_4)=52$. Also in this case $|Z|=1$. Indeed, the existence of an antipodal
map $x \mapsto x'$, would imply the existence of an
antipodal map for every connected totally geodesic submanifold. 
Since $\R P^2$, $\C P^2$ and $\HH P^2$ occur as totally
geodesic subspaces (\cite{W}, Lemma 4), this is not the case.

\medskip

Our next example concerns compact Lie groups viewed as symmetric spaces.

\subsection{Example: Security for compact semisimple Lie groups}  \label{exam2_sub}
Let $K$ be a compact, connected, semisimple Lie group.
Let $<\cdot,\cdot>_k$ be a riemannian metric on $K$ which is both
left-invariant and right-invariant. We will call these metrics
double-invariant. A double-invariant metric is determined by its restriction, 
$<\cdot,\cdot>_e$, to the
tangent space $T_e(K)\sim\k$, which is a $\adj(K)$-invariant inner product.
Conversely, any $\adj(K)$-invariant inner product 
$<\cdot,\cdot>$ on $\k$ uniquely extends to a double-invariant 
riemannian metric on $K$. Although such inner product is not unique,
in general,\footnote{
It is unique, up to scaling, if $K$ is a simple group.}
the results below are valid for any double-invariant metric. In the standard
example, $<\cdot,\cdot>$ is the negative of the Cartan-Killing form. 

We will regard $K$ both as a group and as a symmetric space.
To avoid confusion, we will denote the latter by $[K]$.
The group $K\times K$ acts on $[K]$ via $(k_1,k_2)\cdot[k]=[k_1\,k\,k_2^{-1}]$.
The isotropy group of $[e]$ is the diagonal subgroup
$\{(k,k):\,k\in K\}\subset K\times K$, and we denote it by $K$ as well.
Thus, $[K]=(K\times K)/K$, and $o=[e]$. Let $Z\subset K$ be the center of the group,
and let $Z([K],[e])$ be the center of the symmetric space. 
Then $Z([K],[e])=[Z]$. The mapping $A\mapsto[A]$ provides a one-to-one
correspondence between the set of Cartan subgroups of $K$ and the set
$\ttt_{[e]}([K])$ of maximal tori in $[K]$ containing the reference point.
Hence $\rk(K)=\rk([K])$.

Set $\rk=\rk(K)=\rk([K])$. An element $k\in K$ is regular (resp. singular)
if $\dim(\ker(\adj(k)))=\rk$ (resp. $\dim(\ker(\adj(k)))>\rk$).
By discussion above, a point $[k]$ is regular (resp. maximally singular)
in the sense of Definition~\ref{reg_def}
iff the group element $k$ is regular (resp. $k\in Z$).
Let $\{[k_1],[k_2]\}$ be an arbitrary configuration in $[K]$. In view of the
preceding remarks, $\{[k_1],[k_2]\}$ is a regular configuration
(resp. a maximally singular configuration) iff $k_1k_2^{-1}\in K$
is a regular element (resp. $k_2=z\,k_1$, where $z\in Z$.)

The preceding discussion and Theorem~\ref{block_thm} yield the
following corollary.

\begin{thm}  \label{secur_group_thm}
Let $K$ be a compact, connected, semisimple Lie group
endowed with a double-invariant riemannian metric.
Then the following holds:

\noindent 1. Any configuration $\{k_1,k_2\}$ such that $k_1k_2^{-1}$ is
a regular element of $K$ is secure; it suffices $2^{\rk(K)}$  
points to block all connecting geodesics; 

\noindent 2. There exists an element $z\in Z$ of the centre of $K$ 
such that any configuration $\{k,z\,k\}$ is insecure.
\end{thm}
\section{General compact, locally symmetric spaces}  \label{gener_loc_symm}
We begin with a general proposition. It is of interest by itself;
we will also use it in the proof of Theorem~\ref{gen_sec_thm}.
\begin{prop} \label{loc_iso_prop}
Let $X$ (resp. $Y$) be a (resp. compact) riemannian manifold,
and let $f:Y\to X$ be a local isometry. (We do not assume that
$f$ is onto.) If the space $Y$ is insecure, then $X$ is insecure as well.
\end{prop}
\begin{proof} Let $\{y_1,y_2\}\subset Y$ be an insecure configuration.
We will show that $\{x_1=f(y_1),x_2=f(y_2)\}\subset X$ is also an insecure configuration.
The mapping $f$ sends geodesics in $Y$ into geodesics in $X$.
Let $\Ga_f(x_1,x_2)\subset\Ga(x_1,x_2)$ be the set of connecting geodesics
of the form $\ga=f(\tilde{\ga})$, where $\tilde{\ga}\in\Ga(y_1,y_2)$. 
Suppose that $\{x_1,x_2\}$ is a secure configuration,
and let $F\subset X$ be a blocking set. In particular, $F$ blocks the geodesics
in $\Ga_f(x_1,x_2)$. Set $\tF=f^{-1}(F)$. Then $\tF\subset Y$ is a finite set,
and it blocks all geodesics in $\Ga(y_1,y_2)$, contrary to our assumption.
\end{proof}

Let $X,Y$ be riemannian manifolds. We denote by $X\times Y$ the product
manifold endowed with the product metric. To be precise,
if $(x,y)\in X\times Y$,
then $T_{(x,y)}X\times Y=T_xX\oplus T_yY$. 

\begin{corol}  \label{prod_secur_cor} Let $X,Y$ be arbitrary riemannian manifolds.

\noindent 1. If the space $X\times Y$ is secure, then both $X,Y$ are secure.

\noindent 2. Suppose that $Y\subset X$ is a totally geodesic submanifold.
If $Y$ is insecure then $X$ also is.
\end{corol}
\begin{proof} The first claim follows from the second; the argument of
Proposition~\ref{loc_iso_prop} proofs the second claim. The 
compactness of $Y$ was used in the proof of Proposition~\ref{loc_iso_prop}
only to insure that $|\tF|<\infty$. In our setting $\tF=F$, hence
it is a finite set.
\end{proof}

\medskip

We will now turn to the security of general locally symmetric spaces.
In order to avoid confusion, we will modify our notation for the sets
of connecting geodesics. Namely, we denote by $G(x,y)$ etc the set
of connecting geodesics for the configuration $\{x,y\}$.

\begin{thm}  \label{gen_sec_thm}
Let $M$ be a compact, locally symmetric space. Then $M$ is secure iff
it is of euclidean type.
\end{thm}
\begin{proof} Let $M=\Ga\setminus S$, where $S=\tM$ 
and $\Ga\subset\iso(S)$ is the group of deck transformations. Let 
$p:S\to M$ be the covering map.

The general simply connected, symmetric space has a unique
decomposition $S=S_0\times S_-\times S_+$ \cite{He1},
where some of the factors may be trivial. If $S=S_0\times S_-$, e. g.,
we will say that the factor $S_+$ is not present.
In view of Proposition~\ref{euclid_case}, it suffices to
show that if $S_-$ or $S_+$  are present, then $M$ is insecure.

We will first show that the presence of 
$S_+$ implies the insecurity of $M$. 
Assume the opposite. Then $S=S_+\times Z$. Let $z\in Z$ be arbitrary.
The restriction of $p:S\to M$ to $S_+\times\{z\}\subset S$ satisfies the
assumptions of Proposition~\ref{loc_iso_prop}. Hence, by Theorem~\ref{block_thm} 
and Proposition~\ref{loc_iso_prop}, $M$ is insecure.
 
In view of Theorem~\ref{nct_thm}, it remains to show that if 
$S= S_0 \times S_-$, then $M$ is insecure. In fact, we will prove
that any configuration in $M$ is insecure.
For notational convenience, we set $S_-=S_1$.
Note that $S=S_0\times S_1$ is a Hadamard manifold, and that $S_0$ is
the euclidean de Rham factor.

For $i=0,1$ let
$q_i:S\to S_i$ (resp. $\rho_i: \Gamma \to ISO(S_i)$)
be the natural projections, and set $\Ga_i=\rho_i(\Gamma)$.
By the results of  P. Eberlein \cite{Eb4}, there exists a finite index subgroup
$\Gamma' \subset \Gamma$, such that  
$\Ga_1'=\rho_1(\Gamma')$ is a dicrete, fixed point free, cocompact group 
of isometries  of $S_1$.

Set $M'=\Gamma' \setminus S$.
The inclusion $\Ga'\subset\Ga$ yields a finite covering $\pi:M'\to M$.
By Proposition~\ref{cov_sec_prop}, it suffices to show that
every configuration in $M'$ is insecure. Hence, we assume 
from now on that the group $\Ga$ itself satisfies the conditions above,
and suppress the ``prime" from our notation.

Let $(x_0,x_1)\ne(y_0,y_1)\in S_0\times S_1$ be arbitrary points.
We will use the notation of section~\ref{notcomp_set} for the unique
geodesic, connecting a pair of points of a Hadamard manifold.
The geodesics $[x_0,y_0]\subset S_0,\, [x_1,y_1]\subset S_1$, and
$[(x_0,x_1),(y_0,y_1)]\subset S_0\times S_1$ satisfy
$q_0([(x_0,x_1),(y_0,y_1)])=[x_0,y_0],\,q_1([(x_0,x_1),(y_0,y_1)])=[x_1,y_1]$.
There exist  (unique) linear parametrizations
$z_0(t),z_1(t)$ of the geodesics $[x_0,y_0],[x_1,y_1]$ respectively, such that
$(z_0(t),z_1(t))$ is the arclength parametrization of  the geodesic
$[(x_0,x_1),(y_0,y_1)]$.

Let $\{p(x_0,x_1), p(y_0,y_1)\}\subset M$ be an arbitrary configuration.
Suppose that it is secure, and let $F\subset M$ be a blocking set.
Then for every $\ga',\ga''\in\Ga$ the geodesic $[\ga'(x_0,x_1), \ga''(y_0,y_1)]\subset S$
passes through a point of $\tF=p^{-1}(F)$. 
Let $\ga_i'=\rho_i(\ga'),\ga_i''=\rho_i(\ga'')$, where $i=0,1$.
Then, by preceding remarks, every geodesic $[\ga_1'(x_1),\ga_1''(y_1)]\subset S_1$
passes through a point of $q_1(\tF)$. 

Let $M_1=\Ga_1\setminus S_1$, and let $p_1:S_1\to M_1$ be the covering map.
The set $\tF\subset S$ is a finite union of $\Ga$-orbits.
The projection $q_1:S\to S_1$ sends $\Ga$-orbits onto $\Ga_1$-orbits.
Thus, $q_1(\tF)\subset S_1$ is a finite union of $\Ga_1$-orbits.
Hence, $q_1(\tF)=p_1^{-1}(F_1)$, where $F_1\subset S_1$ is a finite set.

Let $\tG$ be the set of geodesics in $S_1$ given by
$\tG=\{[\ga_1'(x_1),\ga_1''(y_1)]:\ga_1',\ga_1''\in\Ga_1\}$. Then
$\tG=p_1^{-1}(G(p_1(x_1),p_1(y_1)))$. We have observed that
every geodesic in $\tG$ passes through
a point of the set $q_1(\tF)=p_1^{-1}(F_1)$. Since $|F_1|<\infty$,
the set $F_1\subset M_1$ is a blocking set for the configuration $\{(p_1(x_1),p_1(y_1)\}$.
But this contradicts to Theorem~\ref{nct_thm}.
\end{proof}

\medskip

\medskip

We conclude by another characterization of (uniformly) secure, compact,
locally symmetric spaces. It is immediate from Theorem~\ref{gen_sec_thm}
and Proposition~\ref{euclid_case}.

\begin{corol}  \label{gen_sec_cor}
Let $M^n$ be a compact, locally symmetric space. Then the following
claims are equivalent:

\noindent 1. The space $M^n$ is secure; 

\noindent 2. The space $M^n$ is uniformly secure. Its security threshold is bounded 
by a constant that depends only on $n$;

\noindent 3. The space $M^n$ is covered by a flat torus.
\end{corol}

\end{document}